\documentclass{conm-p-l}

\copyrightinfo{2005}{}

\theoremstyle{definition}

\theoremstyle{remark}

\numberwithin{equation}{section}

\newcommand{\tu}{\tilde{u}}

\newcommand{\e}{\epsilon}

\newcommand{\al}{\alpha}
\newcommand{\be}{\beta}

\newcommand{\om}{\omega}
\newcommand{\Om}{\Omega}

\newcommand{\lag}{\langle}
\newcommand{\rag}{\rangle}

\newcommand{\tp}{\tilde{p}}

\begin{document}

\title[On the True Nature of Turbulence]{On the True Nature of Turbulence}

\author{Y. Charles Li}
\address{Department of Mathematics, University of Missouri, 
Columbia, MO 65211}
\curraddr{}
\email{cli@math.missouri.edu}
\thanks{}


\subjclass{Primary 76-02, 37-02, 35-02; Secondary 93-02}



\keywords{Turbulence, Navier-Stokes equations, global well-posedness,
chaos in partial differential equations, control of turbulence}

\begin{abstract}
In this article, I would like to express some of my views on the nature of 
turbulence. These views are mainly drawn from the author's recent results on 
chaos in partial differential equations \cite{Li04}.

Fluid dynamicists believe that Navier-Stokes equations accurately describe
turbulence. A mathematical proof on the global regularity of the solutions to 
the Navier-Stokes equations is a very challenging problem. Such a proof or 
disproof does not solve the problem of turbulence. It may help understanding 
turbulence. Turbulence is more of a dynamical system problem. Studies on 
chaos in partial differential equations indicate that turbulence can have 
Bernoulli shift dynamics which results in the wandering of a turbulent 
solution in a fat domain in the phase space. Thus, turbulence can not 
be averaged. The hope is that turbulence can be controlled.
\end{abstract} 

\maketitle

\tableofcontents








\section{The Governing Equations of Turbulence}

It has been overwhelmingly accepted by fluid dynamicists that the 
Navier-Stokes equations are accurate governing equations of turbulence. 
Their delicate experimental measurements on turbulence have led them to 
such a conclusion. A simple form of the Navier-Stokes equations, 
describing viscous incompressible fluids, can be written as 
\begin{equation}
u_{i,t} + u_j u_{i,j} = - p_{,i} + \ \mbox{Re}^{-1} \ u_{i,jj} +f_i \ , 
\quad u_{i,i} = 0 \ ;
\label{NS}
\end{equation}
where $u_i$'s are the velocity components, $p$ is the pressure, $f_i$'s 
are the external force components, and Re is the Reynolds number. There are 
two ways of deriving the Navier-Stokes equations: (1). The fluid dynamicist's 
way of using the concept of fluid particle and material derivative, (2). 
The theoretical physicist's way of starting from Boltzman equation. 
According to 
either approach, one can replace the viscous term $\mbox{Re}^{-1} \ u_{i,jj}$ 
by for example
\begin{equation}
\mbox{Re}^{-1} \ u_{i,jj} + \al u_{i,jjkk} + \cdots \ . 
\label{hvt}
\end{equation}
Here the only principle one can employ is the Einstein covariance principle
which eliminates the possibility of third derivatives for example. 
According to the fluid dynamicist's way, the viscous term 
$\mbox{Re}^{-1} \ u_{i,jj}$ was derived from a principle proposed by Newton 
that the stress is proportional to the velocity's derivatives (strain, not 
velocity). Such fluids are called Newtonian fluids. Of course, there exist 
non-Newtonian fluids like volcanic lava for which the viscous term is more 
complicated and can be nonlinear. According to the theoretical physicist's way,
the viscous term was obtained from an expansion which has no reason to stop 
at its leading order term $\mbox{Re}^{-1} \ u_{i,jj}$. 

\section{Global Well-Posedness of the Navier-Stokes Equations}

It is well known that the global well-posedness of the Navier-Stokes equations
(\ref{NS}) has been selected by the Clay Mathematics Institute as one of its 
seven one million dollars problems. Specifically, the difficulty lies at the 
global regularity \cite{Ler34}. More precisely, the fact that
\[
\int \int u_{i,j} u_{i,j} \ dx\ dt  
\]
being bounded only implies 
\[
\int u_{i,j} u_{i,j} \ dx
\]
being bounded for almost all $t$, is the key of the difficulty. In fact, 
Leray was able to show that the possible exceptional set of $t$ is 
actually a compact set of measure zero. There have been a lot of more 
recent works on describing this exceptional compact set \cite{CKN82}. 
The claim that this possible exceptional compact set is actually empty, 
will imply the global regularity and the solution of the problem. The 
hope for such a claim seems slim.

Even for ordinary differential equations, often one can not prove their 
global well-posedness, but their solutions on computers look perfectly 
globally regular and sometimes chaotic. Chaos and global regularity are 
compatible. The fact that fluid experimentalists quickly discovered 
shocks in compressible fluids and never found any finite time blow up in 
incompressible fluids, indicates that there might be no finite time 
blow up in Navier-Stokes equations (even Euler equations). On the other 
hand, the solutions of Navier-Stokes equations can definitely be 
turbulent (chaotic). 

Replacing the viscous term $\mbox{Re}^{-1} \ u_{i,jj}$ by higher order 
derivatives (\ref{hvt}), one can prove the global regularity \cite{KP02}.
This leaves the global regularity of (\ref{NS}) a more challenging and 
interesting mathematical problem. Assume that the unthinkable event 
happens, that is, someone proves the existence of a meaningful finite 
time blow up in (\ref{NS}), then fluid experimentalists need to identify 
such a finite time blow up in the experiments. If they fail, then the choice 
will be whether or not to replace the viscous term 
$\mbox{Re}^{-1} \ u_{i,jj}$ in the Navier-Stokes equations (\ref{NS}) 
by higher order derivatives like (\ref{hvt}) to better model the fluid motion.

Even after the global regularity of (\ref{NS}) is proved or disproved, the 
problem of turbulence is not solved although the global regularity 
information will help understanding turbulence. Turbulence is more of a 
dynamical system problem. Often a dynamical system study does not depend on 
global well-posedness. Local well-posedness is often enough. In fact, this 
is the case in my proof on the existence of chaos in partial differential 
equations \cite{Li04}.

\section{Chaos in Partial Differential Equations}

Ever since the discovery of chaos in low dimensional systems, people have 
been trying to use the concept of chaos to understand turbulence 
\cite{RT71}. There are two types of fluid motions: Laminar flows and 
turbulent flows. Laminar flows look regular, and turbulent flows are 
non-laminar and look irregular. Chaos is more precise, for example, in 
terms of Bernoulli shift dynamics. On the other hand, even in low 
dimensional systems, there are solutions which look irregular for a 
while, and then look regular again. Such a dynamics is often called a 
transient chaos.

Everyone knows that the signature of chaos is sensitive dependence on 
initial data. Often the word ``sensitive'' is over-imagined. For any 
fixed large time, the chaotic solution still depends on its initial 
condition continuously. It is the infinite time that leads to sensitive 
dependence. 

Low dimensional chaos is the starting point of a long journey toward 
understanding turbulence. To have a better connection between chaos 
and turbulence, one has to study chaos in partial differential 
equations \cite{Li04}. Take the simple perturbed sine-Gordon equation 
for example \cite{Li04c} \cite{Li05c}
\begin{equation}
u_{tt} = c^2 u_{xx} +\sin u +\e \left [ -a u_t +\cos t \ \sin^3 u \right ]\ ,
\label{PSG}
\end{equation}
which is subject to periodic boundary condition
\[
u(t, x+2\pi ) = u(t,x) \ ,
\]
and even or odd constraint
\[
u(t,-x) =u(t,x) \quad \mbox{or} \quad u(t,-x) =-u(t,x)\ ,
\]
where $u$ is a real-valued function of two real variables ($t,x$), $c$ is 
a real constant, $\e \geq 0$ is a small perturbation parameter, and $a>0$ 
is an external parameter. One can view (\ref{PSG}) as a flow defined in 
the phase space 
\[
(u,u_t) \in H^1 \times L^2 
\]
where $H^1$ and $L^2$ are the Sobolev spaces on [$0,2\pi$]. A point in the 
phase space corresponds to two profiles 
\[
(u(x), u_t(x))\ .
\]
One can prove that there exists a homoclinic orbit $(u,u_t) = h(t,x)$ 
asymptotic to $(u,u_t) = (0,0)$ \cite{Li04c} \cite{Li05c}. Let us define
two orbits segments 
\[
\eta_0 : \ (u,u_t) = (0,0) \ , \quad t \in [-T, T]\ , \quad 
\eta_1 : \ (u,u_t) = h(t,x) \ , \quad t \in [-T, T]\ .
\]
When $T$ is large enough, $\eta_1$ is almost the entire homoclinic orbit 
(chopped off in a small neighborhood of $(u,u_t) = (0,0)$). To any binary 
sequence
\begin{equation}
a = \{ \cdots a_{-2} a_{-1} a_0, a_1 a_2 \cdots \} \ ,
\quad a_k \in \{ 0 , 1 \} \ ; 
\label{bsq}
\end{equation}
one can associate a pseudo-orbit 
\[
\eta_a = \{ \cdots \eta_{a_{-2}} \eta_{a_{-1}} \eta_{a_0}, \eta_{a_1} 
\eta_{a_2} \cdots \} \ .
\]
The pseudo-orbit $\eta_a$ is not an orbit but almost an orbit. One can prove 
that for any such pseudo-orbit $\eta_a$, there is a unique true orbit in 
its neighborhood \cite{Li04c} \cite{Li05c}. Therefore, each binary sequence 
labels a true orbit. All these true orbits together form a chaos. In order to 
talk about sensitive dependence on initial data, one can introduce the 
product topology by defining the neighborhood basis of a binary sequence 
\[
a^* = \{ \cdots a^*_{-2} a^*_{-1} a^*_0, a^*_1 a^*_2 \cdots \}
\]
as 
\[
\Om_N = \left \{ a \ : \quad a_n = a^*_n \ , \quad |n| \leq N \right \}\ .
\]
The Bernoulli shift on the binary sequence (\ref{bsq}) moves the comma 
one step to the right. Two binary sequences in the neighborhood $\Om_N$ 
will be of order $\Om_1$ away after $N$ iterations of the Bernoulli shift. 
Since the binary sequences label the orbits, the orbits will exhibit the 
same feature. In fact, the Bernoulli shift is topologically conjugate 
to the perturbed sine-Gordon flow. 

Replacing a homoclinic orbit by its fattened version -- a homoclinic tube, 
or by a heteroclinic cycle, or by a heteroclinically tubular cycle; one 
can still obtain the same Bernoulli shift dynamics \cite{Li03a} \cite{Li03b} 
\cite{Li04c} \cite{Li05c}. 

Adding diffusive perturbation $\e b u_{txx}$ to (\ref{PSG}), one can still 
prove the existence of homoclinics or heteroclinics, but the Bernoulli 
shift result has not been established \cite{Li04c} \cite{Li05c}.

Another system studied is the complex Ginzburg-Landau equation 
\cite{Li04a} \cite{Li04b},
\begin{equation}
iq_t =q_{xx} +2 \left [ |q|^2 -\om^2 \right ] +
i\e \left [ q_{xx} -\al q +\be \right ] \ ,
\label{CGL}
\end{equation}
which is subject to periodic boundary condition and even constraint 
\[
q(t,x +2\pi ) = q(t,x) \ , \quad q(t,-x)=q(t,x) \ ,
\]
where $q$ is a complex-valued function of two real variables ($t,x$), 
($\om , \al , \be $) are positive constants, and $\e \geq 0$ is a small 
perturbation parameter. In this case, one can prove the existence of 
homoclinic orbits \cite{Li04a}. But the Bernoulli shift dynamics was 
established under generic assumptions \cite{Li04b}.

A real fluid example is the amplitude equation of Faraday water wave, 
which is also a complex Ginzburg-Landau equation \cite{Li04d},
\begin{equation}
iq_t =q_{xx} +2 \left [ |q|^2 -\om^2 \right ] +
i\e \left [ q_{xx} -\al q +\be \bar{q} \right ] \ ,
\label{CGL1}
\end{equation}
subject to the same boundary conditon as (\ref{CGL}). For the first time, 
one can prove the existence of homoclinic orbits for a water wave 
equation (\ref{CGL1}) \cite{Li04d}. The Bernoulli shift dynamics was 
also established 
under generic assumptions \cite{Li04d}. That is, for the first time, one 
can prove the existence of chaos in water waves under generic assumptions.

The nature of the complex Ginzburg-Landau equation is a parabolic equation 
which is near a hyperbolic equation. The same is true for the perturbed 
sine-Gordon equation with the diffusive term $\e b u_{txx}$ added. They 
contain effects of diffusion, dispersion, and nonlinearity. The 
Navier-Stokes equations are diffusion-advection equations. The advective 
term is missing from the perturbed sine-Gordon equation and the complex 
Ginzburg-Landau equation. But the modified KdV equation does contain an 
advective term. In principle, perturbed modified KdV equation should have 
the same feature as the perturbed sine-Gordon equation. Turbulence happens 
when the diffusion is weak, i.e. in the near hyperbolic regime. One should 
hope that turbulence should share some of the features of chaos in the 
perturbed sine-Gordon equation. There is a popular myth that turbulence 
is fundamentally different from chaos because turbulence contains many 
unstable modes. In both the perturbed sine-Gordon equation and the complex 
Ginzburg-Landau equation, one can incorporate as many unstable modes as 
one likes, the resulting Bernoulli shift dynamics is still the same. On a 
computer, the solution with more unstable modes may look rougher, but it
is still chaos. So I think the issue of number of unstable modes between 
turbulence and chaos is an illusion.

Turbulence is any flow that is non-laminar. Sometimes, turbulence can happen 
in a localized spot of a fluid domain, or during a finite period of time. 
These are not chaos. I have a favorite similie of the situation: One can 
think turbulence as marbles; and those flows for which the existence of 
chaos can be rigorously proved, as diamonds. Marbles are everywhere, while 
diamonds are rare. Understanding diamonds can help understanding marbles. 
Diamonds are precious, while marbles are realistically useful in engineering.

A simple setup for studying the chaotic nature of turbulence is posing the 
Navier-Stokes equation (\ref{NS}) on a spatially periodic domain, with a 
temporally and spatially periodic external force. In this case, one can 
take the advantage of Fourier series. One can show that there are 
well-defined invariant manifolds \cite{Li05}. A thorough numerical 
investigation of this dynamical system should be significant for a 
better understanding of turbulence. 

\section{Control of Turbulence}

When dealing with random solutions to a stochastic equation, researchers 
are not content with the random solutions as they are. Various averagings 
will be conducted to gain more certain quantifications of the random 
solutions, since uncertainty is never the favorite to researchers in 
contrast to certainty. Fundamentally encouraging to such thoughts is 
that these averagings are very successful in describing the random solutions.

When dealing with Navier-Stokes equations which are nonlinear deterministic 
equations, fluid engineers are very happy with laminar solutions as they 
are, but not turbulent solutions. They have been trying hard to quantify 
turbulent solutions with averaging techniques. Reynolds envisioned a 
relatively long time averaging to the turbulent solutions. Such an 
averaging failed miserably.

From what we learn about chaos in partial differential equations, turbulent 
solutions not only have sensitive dependences on initial conditions, but 
also are densely packed inside a domain in the phase space. They are far 
away from the feature of fluctuations around a mean. In fact, they wander 
around in a fat domain rather than a thin domain in the phase space. 
Therfore, averaging makes no sense at all. One has to be content with 
turbulent solutions as they are.

In real life, turbulence often represents unpleasant or disastrous events. 
When an airplane meets turbulence, the passengers do not feel comfortable 
and the airplane can be damaged. The fundamental question here is whether 
or not turbulence can be controlled. Here the word ``control'' represents 
a wide spectrum of actions: Taming turbulent states into laminar states 
\cite{ABT01}, reducing turbulent drag \cite{LB98} \cite{Kim03}, enhancing 
turbulent mixing \cite{LB98} \cite{Kim03}, and gearing a turbulent orbit 
to a specific target \cite{CY03} etc.. The final motto that I am aiming 
at is:
\begin{itemize}
\item Turbulence can not be averaged, but can be controlled.
\end{itemize}
Specific control tools have been developed. These are sensors and actuators 
placed in flow fields. These sensors and actuators hopefully can be placed 
by MEMS (Micro-Electro-Mechanical-System) technology in the future to 
obtain a more effective control. 

One can re-interpret the Reynolds averaging as a control of taming 
turbulence into a laminar flow. According to Reynolds, one splits the 
variables in (\ref{NS}) into two parts:
\[
u_i = U_i + \tu_i \ , \quad p=P+\tp
\]
where the capital letters represent relatively long time averages which 
are still a function of time and space, and the tilde-variables represent 
mean zero fluctuations,
\[
U_i = \lag u_i \rag \ , \quad \lag \tu_i \rag =0 \ , \quad
P = \lag p \rag \ , \quad \lag \tp \rag =0 \ .
\]
A better interpretation is by using ensemble average of repeated 
experiments. One can derive the Reynolds equations for the averages,
\begin{equation}
U_{i,t} + U_jU_{i,j} = - P_{,i} + \ \mbox{Re}^{-1} \ U_{i,jj} 
-\lag \tu_i\tu_j \rag_{,j} +f_i \ , 
\quad U_{i,i} = 0 \ .
\label{RE}
\end{equation}
The term $\lag \tu_i\tu_j \rag$ is completely unknown. Fluid engineers 
call it Reynolds stress. The Reynolds model is given by
\begin{equation}
\lag \tu_i\tu_j \rag = -R^{-1} \ U_{i,j} \ ,
\label{RM}
\end{equation}
where $R$ is a constant. There are many more models on the term 
$\lag \tu_i\tu_j \rag$. But no one leads to a satisfactory result. One 
can re-interpret the Reynolds equations (\ref{RE}) as control equations 
of the orginal Navier-Stokes equations (\ref{NS}), with the term 
$\lag \tu_i\tu_j \rag_{,j}$ being the control of taming a turbulent solution 
to a laminar solution (hopefully nearby). The Reynolds model (\ref{RM}) 
amounts to changing the fluid viscosity which can bring a turbulent flow 
to a laminar flow. This laminar flow may not be anywhere near the turbulent 
flow though. Thus, the Reynolds model may not produce satisfactory result 
in comparison with the experiments. Fluid engineers gradually gave up all 
these Reynolds' type models and started directly computing the original 
Navier-Stokes equations (\ref{NS})

An advantage of the control theory is that it can be conducted in a 
trial-correction manner without a detailed knowledge of turbulence. 
Of course, better knowledge of turbulence will help the control. In a 
sense, locating chaos and controlling chaos are intertwined. The Melnikov 
intergal can predict the existence of chaos \cite{Li04}, at the same time, 
it also predicts the non-existence of chaos when parameters are changed.

\end{document}